\documentclass{amsart}
\usepackage{amscd,amssymb}
\usepackage[all]{xy}

\newcommand{\rad}{\operatorname{rad}\nolimits}

\newcommand{\Hom}{\operatorname{Hom}\nolimits}

\renewcommand{\Im}{\operatorname{Im}\nolimits}
\newcommand{\Ker}{\operatorname{Ker}\nolimits}

\newcommand{\rrad}{\mathfrak{r}}

\newcommand{\Ext}{\operatorname{Ext}\nolimits}

\newcommand{\op}{{\operatorname{op}\nolimits}}

\newcommand{\frako}{\mathfrak{o}}
\newcommand{\frakt}{\mathfrak{t}}

\renewcommand{\L}{\Lambda}

\newcommand{\B}{{\mathcal B}}

\newcommand{\I}{{\mathcal I}}

\newcommand{\extto}{\xrightarrow}

\newtheorem{lem}{Lemma}[section]
\newtheorem{prop}[lem]{Proposition}
\newtheorem{cor}[lem]{Corollary}
\newtheorem{thm}[lem]{Theorem}
\newtheorem*{Theorem}{Theorem}
\theoremstyle{definition}

\begin{document}

\title{Resolutions over Koszul algebras}
\author[Green]{Edward L. Green}
\thanks{The three first authors received financial support from a
  joint project of CNPq-NSF. The last author thanks the Fapesp-Brasil
  for financial support for a scientific visit to Brasil where this
  work was done.  Finally the third authour thanks CNPq for a research
  grant.}  
\address{Edward L. Green\\ Department of Mathematics\\
Virginia Tech\\ Blacksburg, VA 24061\\ USA}
\email{green@math.vt.edu}
\author[Hartman]{Gregory Hartman}
\address{Gregory Hartman\\ Department of Mathematics\\
The University of Arizona\\ 
617 N. Santa Rita Ave.\\
P.O. Box 210089\\
Tucson, AZ 85721-0089\\ USA}
\email{hartman@math.arizona.edu}
\author[Marcos]{Eduardo N. Marcos}
\address{Eduardo N. Marcos\\ 
Instituto de Matem\'atica e Estat\'istica\\
Universidade S\~ao Paulo (IME-USP)\\
Rua do Mat\~ao, 1010 - Cidade Universit\'aria\\
CEP 05508-090\\
S\~ao Paulo - SP - Brazil}
\email{enmarcos@ime.usp.br}
\author[Solberg]{\O yvind Solberg}
\address{\O yvind Solberg\\Institutt for matematiske fag\\
NTNU\\ N--7491 Trondheim\\ Norway}
\email{oyvinso@math.ntnu.no}

\date{\today}

\begin{abstract}
In this paper we show that if $\L=\amalg_{i\geq 0}\L_i$ is a Koszul
algebra with $\L_0$ isomorphic to a product of copies of a field, then
the minimal projective resolution of $\L_0$ as a right $\L$-module
provides all the information necessary to construct both a minimal
projective resolution of $\L_0$ as a left $\L$-module and a minimal
projective resolution of $\L$ as a right module over the enveloping
algebra of $\L$. The main tool for this is showing that there is a
comultiplicative structure on a minimal projective resolution of
$\L_0$ as a right $\L$-module.
\end{abstract}

\maketitle

\section*{Introduction and preliminaries}
Let $\L=\amalg_{i\geq 0} \L_i$ be a Koszul algebra over a field $k$
with $\L_0$ a product of copies of $k$, where we recall the definition
of Koszul later in this section.  Denote by $(\mathbb{L},e)$ a minimal
(graded) projective resolution of $\L_0$ as a right $\L$-module.  We
show that $(\mathbb{L},e)$ contains all the information needed to
construct a minimal projective resolution of $\L$ as a right
$\L^e$-module, where $\L^e=\L^\op\otimes_k \L$. The resolution
$(\mathbb{L},e)$ is shown to have a ``comultiplicative
structure''. This structure is used to prove that one can obtain a
minimal projective resolution of $\L_0$ over $\L$ as a left
$\L$-module from the knowledge of $(\mathbb{L},e)$.  We apply these
results to prove an unpublished result of E. L. Green and D. Zacharia
that $\L$ is a Koszul algebra if and only if $\L$ is a linear module
as a right module over $\L^e$. In \cite{BGS}, the comultiplicative
structure is applied to give the multiplicative structure of the
Hochschild cohomology ring of a Koszul algebra and also the structure
constants for a basis for the Koszul dual.

The rest of the section is devoted to recalling definitions, results,
and terminology relevant to this paper.  Let $\L=\amalg_{i\geq 0}
\L_i$ be a graded algebra over a field $k$. Assume that (i) $\L_0$ is
a product of copies of $k$, that (ii) each $\L_i$ is finite
dimensional over $k$, and that (iii) $\L$ as an algebra is generated in
degrees $0$ and $1$. Such an algebra $\L$ is isomorphic to a quotient
of the path algebra $kQ/I$, where $kQ$ is isomorphic to the tensor
algebra $T_{\L_0}(\L_1)=\amalg_{i\geq
  0}\underbrace{\L_1\otimes_{\L_0}\cdots \otimes_{\L_0}\L_1}_i$.
Conversely, if $Q$ is a quiver and $I$ is an ideal generated by length
homogeneous elements in $kQ$, then $\L=kQ/I$ is a graded algebra over
$k$ satisfying the conditions above. Throughout this paper $\L$ denotes
a graded algebra having properties (i)--(iii).

Let $\rrad =\amalg_{i\geq 1}\L_i$, which is the graded Jacobson
radical of $\L$.  If $(\mathbb{P},d)$: 
\[\cdots\to P^2\extto{d^2} P^1\extto{d^1} P^0\extto{d^0} M\to 0\]
is a graded projective resolution of a graded $\L$-module $M$, then it
is \emph{minimal} if $\Im d^n\subseteq \rrad P^{n-1}$ for $n\geq
1$. It is well known that graded modules over graded algebras have
minimal graded projective resolutions. We say that a graded projective
resolution 
\[\cdots\to P^2\extto{d^2} P^1\extto{d^1} P^0\extto{d^0} M\to 0\]
is \emph{linear}, and $M$ is a \emph{linear module} if, for $n\geq 0$,
the graded module $P^n$ is generated in degree $n$.  Note that a
linear resolution is a minimal projective resolution.  A graded
algebra $\L$ is a \emph{Koszul algebra} if $\L_0$ is a linear 
module; that is, $\L_0$ has a linear (graded) projective
resolution $(\mathbb{L},e)$:
\[\cdots\to L^2\extto{e^2} L^1\extto{e^1} L^0\extto{e^0} \L_0\to 0\]
as a right $\L$-module.

Before giving the precise results, we introduce notation and recall
results from \cite{GSZ} which are used throughout the paper.  For ease
of notation, let $R=kQ$, let $\B$ be the set of all paths in the
quiver $Q$, and denote by $\B_t$ all the paths of length $t$.

There exist integers $\{ t_n\}_{n\geq 0}$ and elements $\{
f^n_i\}_{i=0}^{t_n}$ in $R$ such that a minimal right projective
resolution $(\mathbb{L},e)$ of $\L_0$ can be given in terms of a
filtration of right ideals
\[\cdots\subseteq \amalg_{i=0}^{t_n} f^n_iR \subseteq \amalg_{i=0}^{t_{n-1}} 
f^{n-1}_iR\subseteq \cdots \subseteq \amalg_{i=0}^{t_1}
f^1_iR\subseteq \amalg_{i=0}^{t_0}f^0_iR=R\] in $R$. Then
$L^n=\amalg_{i=0}^{t_n} f^n_iR/\amalg_{i=0}^{t_n} f^n_iI$ and the
differential $e$ is induced by the inclusion $\amalg_{i=0}^{t_n}
f^n_iR \subseteq \amalg_{i=0}^{t_{n-1}} f^{n-1}_iR$. This inclusion
gives elements $h^{n-1,n}_{ji}$ in $R$ such that
\[f^n_i=\sum_{j=0}^{t_{n-1}} f^{n-1}_jh^{n-1,n}_{ji}\]
for all $i=0,1,\ldots,t_n$ and all $n\geq 1$, so that 
\[e^n(\overline{f^n_i})=
(\overline{h^{n-1,n}_{0i}},\overline{h^{n-1,n}_{1i}},\ldots, 
\overline{h^{n-1,n}_{t_{n-1}i}})\]
for all $n\geq 1$, where $\overline{*}$ denotes the natural residue
class of $*$ modulo $I$. It is shown in \cite{GSZ} that the $f^n_i$'s
can be chosen so that $(\mathbb{L},e)$ is a minimal resolution of
$\L_0$ over $\L$. We point out that an algorithmic construction of
the elements $f^n_i$'s can be found in \cite{G}.

An important property of the elements $\{ f^n_i\}_{i=0}^{t_n}$ is that
there exist elements ${f^{n+1}_j}'$ in $\amalg_{i=0}^{t_{n-1}}f^n_iI$
such that 
\[(\amalg_{i=0}^{t_n} f^n_iR)\cap (\amalg_{i=0}^{t_{n-1}}
f^n_iI)=(\amalg_{i=0}^{t_{n+1}} f^{n+1}_iR) \amalg (\amalg_j
{f^{n+1}_j}'R).\] 

Recall that an element $x$ in $R$ is called \emph{uniform} if $x$ is
non-zero and there exist vertices $u$ and $v$ in $Q$ such that
$x=uxv$. If $x$ is a uniform element with $x=uxv$, then we write
$\frako(x)=u$ and $\frakt(x)=v$.  The elements $f^n_i$ can all be
chosen uniform for $i=0,1,\ldots,t_n$ and all $n\geq 0$, and we assume
that they are.

Note that $t_0+1$ is the number of non-isomorphic graded simple right
$\L$-modules, and that $\{ f^0_i\}_{i=0}^{t_0}$ is the set of vertices
of $Q$. Moreover, $t_1+1$ is the number of arrows of $Q$ and $\{
f^1_i\}_{i=0}^{t_1}$ is choosen to be the set of arrows of $Q$. The
set $\{ f^2_i\}_{i=0}^{t_2}$ is a set of uniform length homogeneous
minimal generators for $I$.

In case $\L$ is a Koszul algebra, we have the following additional
property of the elements $f^n_i$ in $R$; namely each $f^n_i$ is a
linear combination of paths in $\B_n$ for $i=0,1,\ldots,t_n$ and the
length of each path occurring in ${f^n_i}'$ is at least $n+1$. By
length considerations, $h^{n-1,n}_{ji}$ are all linear combinations of
elements in $\B_1$.

In section \ref{section:1} we prove that the elements
$\{f^n_i\}_{i=0,n\geq 0}^{t_n}$ have the following ``comultiplicative
structure'', which is used in \cite{BGS} to give the multiplicative
structure of the Hochschild cohomology ring of a Koszul algebra and
the structure constants for the basis associated to the elements
$\{f^n_i\}$ for the Koszul dual.

\begin{Theorem}
Let $\L=kQ/I$ be a Koszul algebra. Then for each $r$, with $0\leq
r\leq n$, and $i$, with $0\leq i\leq t_n$, there exist elements
$c_{pq}(n,i,r)$ in $k$ such that
\[f^n_i  = \sum_{p=0}^{t_r}\sum_{q=0}^{t_{n-r}} c_{pq}(n,i,r)f^r_pf^{n-r}_q\]
for all $n\geq 1$, all $i$ in $\{0,1,\ldots,t_n\}$ and all $r$ in
$\{0,1,\ldots,n\}$. 
\end{Theorem}
Viewing $\L_0$ as a left module over $\L$, it also has a minimal
graded projective resolution given by $\{ g^n_i\}_{i=0}^{s_n}$, where
$g^n_i$'s are the left analogue of the right $f^n_j$'s in $R$. The
above result is used to prove that one can choose the elements
$g^n_i$'s to be the same as the elements $f^n_j$'s and then the
formula $f^n_i=\sum_{p,q=0}^{t_1,t_{n-1}} c_{pq}(n,i,1)f^1_pf^{n-1}_q$
gives the differential in the projective resolution of $\L_0$ as a
left $\L$-module. Thus the knowledge of the minimal projective
resolution $(\mathbb{L},e)$ via the elements $\{ f^n_i\}$ contains all
the information needed to construct a minimal projective resolution of
$\L_0$ as a left $\L$-module.

In the final section of the paper, the elements $f^n_i$'s are shown to
provide all the information needed to construct a minimal projective
resolution of $\L$ as a right $\L^e$-module.  In particular, we prove
the following.
\begin{Theorem}
Let $\L=kQ/I$ be a Koszul algebra, and let $\{f^n_i\}_{i=0}^{t_n}$ be
defined as above for $\L_0$ as a right $\L$-module. A minimal
projective resolution $(\mathbb{P},\delta)$ of $\L$ over $\L^e$ is
given by 
\[P^n=\amalg_{i=0}^{t_n} \L\frako(f^n_i)\otimes_k \frakt(f^n_i)\L\] 
for $n\geq 0$, where $j$-th component of the differential $\delta^n\colon
P^n\to P^{n-1}$ applied to the $i$-th generator
$\frako(f^n_i)\otimes\frakt(f^n_i)$ is given by  
\[\sum_{p=0}^{t_1} c_{pj}(n,i,1)\overline{f^1_{p}}\frako(f^{n-1}_j)
\otimes\frakt(f^{n-1}_j)+(-1)^n \sum_{q=0}^{t_1}
c_{jq}(n,i,n-1)\frako(f^{n-1}_j)\otimes
\frakt(f^{n-1}_j)\overline{f^1_{q}}\]
for $j=0,1,\ldots,t_{n-1}$ and $n\geq 1$, and $\delta^0\colon
\amalg_{i=0}^{t_0}\L e_i\otimes_k e_i\Lambda \to \L$ is the
multiplication map.
\end{Theorem}
As mentioned earlier, the final result of the paper is that $\L$ is a
Koszul algebra if and only if $\L$ is a linear module as a right
$\L^e$-module.

\section{A resolution with comultiplicative structure}\label{section:1}

In this section Theorem \ref{thm:comult} provides a comultiplicative
structure to a minimal projective resolution of $\L_0$ as a right
$\L$-module. This result is then applied to show that the knowledge of
a minimal projective resolution of $\L_0$ as a right $\L$-module
is sufficient to construct a minimal projective resolution of
$\L_0$ as a left $\L$-module.

Let $\L=kQ/I$ be a graded algebra over a field $k$. Let $\{
t_n\}_{n\geq 0}$ and $\{ f^n_i\}_{i=0,n\geq 0}^{t_n}$ be as in the
introduction. We say that \emph{$\{ f^n_i\}_{i=0,n\geq 0}^{t_n}$
defines a minimal resolution} if the resolution described in the
introduction is minimal. 

The next result shows that the elements $\{ f^n_i\}$ have a
comultiplicative structure for a Koszul algebra. 

\begin{thm}\label{thm:comult}
Let $\L=kQ/I$ be a Koszul algebra, and assume that
$\{f^n_i\}_{i=0}^{t_n}$ defines a minimal resolution of $\L_0$ as
a right $\L$-module. Then for each $r$, with $0\leq r\leq n$, and $i$, 
with $0\leq i\leq t_n$, there exist elements $c_{pq}(n,i,r)$ in $k$
such that
\[f^n_i=\sum_{p=0}^{t_r}\sum_{q=0}^{t_{n-r}} c_{pq}(n,i,r)f^r_pf^{n-r}_q.\]
\end{thm}
\begin{proof}
For any $n$, and $r$ equal to $0$ or $n$, the result follows from
$f^n_i=f^n_if^0_{\frakt(f^n_i)}=f^0_{\frako(f^n_i)}f^n_i$ for
$i=0,1,\ldots, t_n$. Also, this proves the result in the case $n$ is 
equal to $1$.

Next we discuss the case $n=2$. As we have remarked, each
$f^2_i=\sum_{j=0}^{t_1} f^1_jh^{1,2}_{ji}$. Since $\L$ is Koszul, each
$f^2_i$ is a linear combination of paths in $\B_2$, and hence
$h^{1,2}_{ji}$ is a linear combination of elements in $\B_1$. This
gives the result for $n=2$.

Now we proceed by induction on $n$ and assume that the result is true
for $l<n$ and $n\geq 3$. We have that $f^n_i=\sum_{j=0}^{t_{n-1}}
f^{n-1}_jh^{n-1,n}_{ji}$. As in our discussion for $n=2$, we see that
$h^{n-1,n}_{ji}$ is a linear combination of elements in $\B_1$. There
exist elements $c_{ijs}$ in $k$ such that
\[f^n_i=\sum_{j=0}^{t_{n-1}}\sum_{s=0}^{t_1} c_{ijs}f^{n-1}_jf^1_s.\]
By induction, there exist elements $c_{juv}'$ in $k$ such that  
\[f^{n-1}_j=\sum_{u=0}^{t_r}\sum_{v=0}^{t_{n-r-1}}
c'_{juv}f^r_uf^{n-r-1}_v\]
for any $r$, with $0\leq r\leq n-1$. Hence 
\[f^n_i=\sum_j\sum_s\sum_u\sum_v c_{ijs}c'_{juv}f^r_uf^{n-r-1}_vf^1_s\]
for any $r$, with $0\leq r\leq n-1$. The term after $f^r_u$ is 
\begin{equation}
A=\sum_j\sum_s\sum_v c_{ijs}c'_{juv}f^{n-r-1}_vf^1_s.\label{eq:1}
\end{equation}
Theory tells us that 
\[f^n_i=\sum_{w=0}^{t_{n-2}} f^{n-2}_w z_w,\]
where $z_w$ is in $I$. Again by length considerations each $z_w$ is a
linear combination of $f^2_l$'s. Hence, there exist elements
$c_{iwx}''$ in $k$ such that 
\[f^n_i=\sum_{w=0}^{t_{n-2}}\sum_{x=0}^{t_2} c''_{iwx}
f^{n-2}_wf^2_x.\]
By induction each $f^{n-2}_w$ is a linear combination of
$f^r_uf^{n-r-2}_y$. We obtain 
\[f^n_i=\sum_w\sum_x\sum_u\sum_y c''_{iwx}c'''_{wuy} f^r_uf^{n-r-2}_yf^2_x\]
for some $c'''_{wuy}$ in $k$. So the term after $f^r_u$ in this
expression is 
\begin{equation}
B=\sum_w\sum_x\sum_u c''_{iwx}c'''_{wuy} f^{n-r-2}_yf^2_x.\label{eq:2}
\end{equation}
Since $\sum_u f^r_uR$ is a direct sum, we see that formulas 
\eqref{eq:1} and \eqref{eq:2} are equal. The equation \eqref{eq:1}
implies that $A$ is in $\amalg_{v=0}^{t_{n-r-1}} f^{n-r-1}_vR$, and
the equation \eqref{eq:2} implies that $A$ is in
$\amalg_{y=0}^{t_{n-r-2}} f^{n-r-2}_yI$.  It follows that $A$ is
contained in $(\amalg_{t=0}^{t_{n-r}} f^{n-r}_t R) \amalg (\amalg_l
{f^{n-r}_l}' R)$. By length arguments we infer that $A$ is in
$\amalg_{t=0}^{t_{n-r}} f^{n-r}_t R$ and that $A$ is a $k$-linear
combination of the $f^{n-r}_t$'s. Hence we conclude that $f^n_i$ is a
$k$-linear combination of $f^r_s f^{n-r}_t$, and this completes the
proof of the result.
\end{proof}

Since the maps in the minimal projective resolution of $\L_0$ as a
right $\L$-module are given by the $h^{n-1,n}_{ji}$, we explicitly
point out the following relationship. 

\begin{cor}
Keeping the notation of \emph{Theorem \ref{thm:comult}}, we have
\[h^{n-1,n}_{ji}=\sum_{l=0}^{t_1} f^1_lc_{jl}(n,i,n-1)\]
for $n\geq 1$, and $i$ and $j$, with $0\leq i\leq t_n$, and $0\leq j\leq
t_{n-1}$.
\end{cor}

Before applying Theorem \ref{thm:comult} we need the following lemma,
where $J$ denotes the ideal generated by the arrows in $Q$. 

\begin{lem}\label{lem:linindepdirect}
Let $\{ x_i\}_{i\in\I}$ be a set of elements in the linear span of
$\B_s$.  Suppose that $\{ x_i\}_{i\in\I}$ is linearly independent
viewed as vectors over $k$. Then $\sum_{i\in\I} Rx_i$ and
$\sum_{i\in\I} x_iR$ are direct sums.
\end{lem}
\begin{proof}
Suppose that $\sum_i\sum_{j\in\I} c_{ij}q_{ij}x_j=0$ in $R$ for some
elements $c_{ij}$ in $k$ and some paths $q_{ij}$ in $R$. Since all the
paths occurring in any $x_i$ have the same length, we can assume
without loss of generality that the paths $q_{ij}$ all have the same
length, say $t$.  Since $\displaystyle{J^t=\amalg_{q\in B_t} qR}$,
fixing $q$ in $\B_t$, it follows that $\sum_{q_{ij}=q}
q(c_{ij}x_j)=0$, which implies that $\sum_{q_{ij}=q} c_{ij}x_j=0$. By
assumption, we have that $c_{ij}=0$ for all $q_{ij}=q$. Hence we infer
that $\sum_{i\in\I} Rx_i$ is a direct sum. Similarly, $\sum_{i\in\I}
x_iR$ is a direct sum.
\end{proof}

We now show that the $\{ f^n_i\}$ obtained from a right minimal projective
resolution of $\L_0$ and the $\{ g^n_j\}$ obtained from a left
minimal projective resolution of $\L_0$ can be chosen to be the
same. 

\begin{prop}\label{prop:leftright}
Let $\L=kQ/I$ be a Koszul algebra. Let $\{f^n_i\}_{i=0}^{t_n}$ and
$\{g^n_i\}_{i=0}^{s_n}$ define a minimal resolution of $\L_0$ as a
right $\L$-module and as a left $\L$-module, respectively. Then
$s_n=t_n$ for all $n\geq 0$ and the set $\{ g^n_i\}_{i=0}^{t_n}$ can
be chosen to be equal to the set $\{ f^n_i\}_{i=0}^{t_n}$ for all
$n\geq 0$.
\end{prop}
\begin{proof}
For $n$ equal to $0$, $1$, or $2$ the result is clear. Let $n\geq
3$. We proceed by induction on $n$ and assume that the result is true
for all $i<n$.  By Theorem \ref{thm:comult}, for each $i$ with $0\leq
i\leq t_n$ the equalities 
\[f^n_i=\sum_{p,q} c_{pq}f^1_pf^{n-1}_q =\sum_{p',q'} c_{p',q'}
f^2_{p'}f^{n-2}_{q'}\] hold for some $c_{pq}$ and $c_{p',q'}$ in $k$.
Hence $f^n_i$ is in $(\amalg_q Rf^{n-1}_q)\cap
(\amalg_{q'}If^{n-2}_{q'})$, which, by induction, is equal to
$(\amalg_{i=0}^{s_n} Rg^n_i)\amalg (\amalg_i R{g^n_i}')$. Since
$\sum_i f^n_iR$ is direct, the set $\{f^n_i\}$ is linearly independent
as vectors over $k$, and therefore the sum $\sum_i Rf^n_i$ is direct
by Lemma \ref{lem:linindepdirect}. By length considerations, $\{
f^n_i\}_{i=0}^{t_n}$ is contained in the $k$-linear span of $\{
g^n_i\}_{i=0}^{s_n}$. Therefore $t_n\leq s_n$. By switching the roles
of $\{ f^n_i\}$ and $\{ g^n_j\}$ and using the argument above, we
conclude that $\{ g^n_i\}_{i=0}^{s_n}$ is linearly independent and
each $g^n_i$ is in $k$-linear span of $\{ f^n_i\}_{i=0}^{t_n}$. Hence
$s_n=t_n$. By Lemma \ref{lem:linindepdirect} it follows that
$\amalg_{i=0}^{t_n} Rf^n_i=\amalg_{i=0}^{t_n} Rg^n_i$. This shows that
we can choose the set $\{ g^n_i\}_{i=0}^{t_n}$ equal to $\{
f^n_i\}_{i=0}^{t_n}$.
\end{proof}

Proposition \ref{prop:leftright} implies that, given a minimal
projective resolution of $\L_0$ as a right $\L$-module in the from
of $\{ f^n_i\}$, we have all the information to construct a minimal
projective resolution of $\L_0$ as a left $\L$-module. More
precisely, take the $\{ f^n_i\}$ as the $\{ g^n_i\}$, and the maps in
the left resolution are given by $g^n_i\mapsto
\sum_{p=0}^{t_1}\sum_{q=0}^{t_{n-1}} c_{pq}(n,i,1)g^1_pg^{n-1}_q$.

\section{A minimal projective bimodule resolution of $\L$} 

In this section we turn our attention to the construction of a minimal
projective $\L^e$-resolution of $\L$. This construction uses the
comultiplicative structure of the minimal projective resolution of
$\L_0$ as a right $\L$-module found in Theorem
\ref{thm:comult}. This is applied to show an unpublished result of
E. L. Green and D. Zacharia that $\L$ is a Koszul algebra if and only
if $\L$ is a (right) linear module over $\L^e$. 

The following result also shows that the knowledge of the $\{ f^n_i\}$
from a minimal projective resolution of $\L_0$ as a right
$\L$-module is sufficient to explicitly give the projective modules
and the differentials in a minimal projective resolution of $\L$ as a
right $\L^e$-module. The structure of the projective modules in a
minimal projective resolution of $\L$ as right $\L^e$-module was first
given in \cite{H}. Recall that the notation $\overline{*}$ denotes
the natural residue class of $*$ modulo $I$.  Let $\{c_{pq}(n,i,r)\}$
be as in Theorem \ref{thm:comult}. 

\begin{thm}\label{thm:bimodres}
Let $\L=kQ/I$ be a Koszul algebra, and let $\{f^n_i\}_{i=0}^{t_n}$
define a minimal resolution of $\L_0$ as a right $\L$-module. A
minimal projective resolution $(\mathbb{P},\delta)$ of $\L$ over
$\L^e$ is given by
\[P^n=\amalg_{i=0}^{t_n} \L\frako(f^n_i)\otimes_k \frakt(f^n_i)\L\] 
for $n\geq 0$, where the $j$-th component
of the differential $\delta^n\colon P^n\to P^{n-1}$ 
applied to the $i$-th generator
$\frako(f^n_i)\otimes\frakt(f^n_i)$ is given by  
\[\sum_{p=0}^{t_1} c_{pj}(n,i,1)\overline{f^1_{p}}\frako(f^{n-1}_j)
\otimes\frakt(f^{n-1}_j) +(-1)^n \sum_{q=0}^{t_1}
c_{jq}(n,i,n-1)\frako(f^{n-1}_j)\otimes
\frakt(f^{n-1}_j)\overline{f^1_{q}}\]
for $j=0,1,\ldots,t_{n-1}$ and $n\geq 1$, and $\delta^0\colon
\amalg_{i=0}^{t_0}\L e_i\otimes_k e_i\Lambda \to \L$ is the
multiplication map.

In particular, $\L$ is a linear module over $\L^e$. 
\end{thm}
\begin{proof}
Direct computations show that $(\delta)^2=0$, so that
$(\mathbb{P},\delta)$ is a linear complex. In addition, note that
$(\L_0\otimes_\L\mathbb{P},1_{\L_0}\otimes\delta)$ is a
minimal resolution of $\L_0$ as a right $\L$-module. 

In our setting, we have that $\L^e/\rad \L^e\simeq
\Hom_k(\L_0,\L_0)$. Let $(\mathbb{F},d)$ be a minimal
projective resolution of $\L$ as a right $\L^e$-module. Then by
\cite[Chap.\ IX, Proposition 4.3]{CE} we have that
\begin{align}
\Hom_{\L^e/\rad\L^e}(F^n/F^n\rad\L^e,\L^e/\rad\L^e) & \simeq 
\Ext^n_{\L^e}(\L,\L^e/\rad\L^e)\notag\\
& \simeq \Ext^n_\L(\L_0,\L_0)\notag\\
& \simeq \Hom_{\L_0}(\amalg_{i=0}^{t_n} f^n_iR/\amalg_{i=0}^{t_n}
f^n_iJ, \L_0)\notag
\end{align}
for all $n\geq 0$. In particular, $P^n\simeq F^n$ as $\L^e$-modules
for all $n\geq 0$, and hence $P^n/P^n\rad\L^e\simeq F^n/F^n\rad\L^e$
as $\L^e/\rad\L^e$-modules for all $n\geq 0$. Note that these need not
be isomorphic as graded modules, but, we in fact show that this is the
case.

Since $(\mathbb{P},\delta)$ is a complex, we obtain the following
commutative diagram
\[\xymatrix{%
\cdots \ar[r] & P^2\ar[r]^{\delta^2} \ar[d]^{\alpha^2} & 
P^1\ar[r]^{\delta^1} \ar[d]^{\alpha^1} & 
P^0\ar[r]^{\delta^0} \ar[d]^{\alpha^0} & \L \ar[r]\ar@{=}[d] & 0\\
\cdots \ar[r] & F^2\ar[r]^{d^2} & 
F^1\ar[r]^{d^1} & 
F^0\ar[r]^{d^0} & \L \ar[r]& 0
}\]
\sloppy Clearly $\alpha^0\colon P^0\to F^0$ is an isomorphism, and we get an
isomorphism $\alpha^0|_{\Ker\delta^0}\colon \Ker\delta^0\to \Ker
d^0$. Hence $\Ker\delta^0/\Ker\delta^0\rad\L^e\simeq
F^1/F^1\rad\L^e$. Since $\Im\delta^1$ is contained in $\Ker \delta^0$,
this induces a map $\beta^1\colon P^1_1\to
\Ker\delta^0/\Ker\delta^0\rad\L^e$.  If $\beta^1$ is an isomorphism,
then $\alpha^1$ is an isomorphism and we have exactness at $P^0$. Suppose that
$\beta^1$ is not an isomorphism. Since $P^1$ is generated in degree
$1$, there is some projective summand of $P^1$ which is mapped to zero
by $\delta^1$. Using the observation that
$(\L_0\otimes_\L\mathbb{P},1_{\L_0}\otimes\delta)$ is a
minimal resolution of $\L_0$ as a right $\L$-module, we obtain a
contradiction. Hence $\beta^1$ is an isomorphism. 

Since $\alpha^1$ is an isomorphism, we can use the above argument
replacing $\alpha^0$ by $\alpha^1$ to show that $\alpha^2$ is an
isomorphism and exactness at $P^1$. By induction we infer that
$(\mathbb{P},\delta)$ is exact. Since the terms
$\overline{f^1_{p}}\frako(f^{n-1}_j) \otimes\frakt(f^{n-1}_j)$ and
$\frako(f^{n-1}_j)\otimes \frakt(f^{n-1}_j)\overline{f^1_{q}}$ are
elements of degree one in $\L^e$, we conclude that
$(\mathbb{P},\delta)$ is a minimal linear projective resolution of
$\L$ over $\L^e$. This also implies that $\L$ is a linear module over
$\L^e$. The proof is now complete.
\end{proof}

As a consequence we obtain the next corollary which was first proved
by E. L. Green and D. Zacharia.

\begin{cor}
Let $\L=kQ/I$ be a graded algebra. Then $\L$ is a Koszul algebra if
and only if $\L$ as a (right) $\L^e$-module is a linear module. 
\end{cor}
\begin{proof}
Suppose that $\L$ is a Koszul algebra. Then Theorem \ref{thm:bimodres}
implies that $\L$ has a linear projective $\L^e$-resolution, and hence
$\L$ is a linear module over $\L^e$.

Suppose that $\L$ is a linear module over $\L^e$ as a right
module. Let $(\mathbb{P},\delta)$ be a linear projective resolution of
$\L$ as a right $\L^e$-module. Tensoring $\mathbb{P}$ with $\L_0$,
we obtain $\L_0\otimes_\L\mathbb{P}$. But
$\L_0\otimes_\L\mathbb{P}$ is a linear projective resolution of
$\L_0$ as a right $\L$-module. Hence $\L$ is a Koszul algebra, and
we are done.
\end{proof}

\end{document}